\documentclass[twoside,12pt]{article}
\usepackage{amssymb,amsmath,bm,euscript}
\usepackage{graphicx,color,caption}
\usepackage{ifpdf}

\usepackage[colorlinks=true]{hyperref}

\usepackage{braket}

\usepackage{algorithm, algorithmic}

\newtheorem{proposition}{Proposition}[section]
\newtheorem{theorem}[proposition]{Theorem}

\newtheorem{corollary}[proposition]{Corollary}

\newcommand{\qed}{\hphantom{.}\hfill $\Box$\medbreak}

\def\S{\mathcal{S}}

\def\S{\mathcal{S}}

\def\0{{\bf 0}}

\title{\bf{Positivity Conditions for Cubic, Quartic and Quintic Polynomials}}%\thanks{This research was supported by the Hong Kong
\author{ \hspace{1mm} Liqun Qi\thanks{%Department of Mathematics, School of Science, Hangzhou Dianzi University, Hangzhou 310018 China; Future Network Theory Lab, 2012 Labs
Department of Applied
    Mathematics, The Hong Kong Polytechnic University, Hung Hom,
    Kowloon, Hong Kong, China; ({\tt liqun.qi@polyu.edu.hk}).},
     \  \
   Yisheng Song\thanks{School of Mathematical Sciences, Chongqing Normal University, Chongqing 401331 China; ({\tt yisheng.song@cqnu.edu.cn}). This author's work was supported by NSFC (Grant No.  11571095, 11601134). },
 \ and \
  Xinzhen Zhang\thanks{School of Mathematics, Tianjin University, Tianjin 300354 China; ({\tt xzzhang@tju.edu.cn}). This author's work was supported by NSFC (Grant No.  11871369). },
}

\begin{document}
\date{\today}
\maketitle

\begin{abstract}
We present a necessary and sufficient condition for a cubic polynomial to be positive for all positive reals.  We identify the set where the cubic polynomial is nonnegative but not all positive for all positive reals, and explicitly give the points where the cubic polynomial attains zero.  We then reformulate a necessary and sufficient condition for a quartic polynomial to be nonnegative for all positive reals.  From this, we derive a necessary and sufficient condition for a quartic polynomial to be nonnegative and positive for all reals.   Our condition explicitly exhibits the scope and role of some coefficients, and has strong geometrical meaning.   In the interior of the nonnegativity region for all reals, there is an appendix curve.   The discriminant is zero at the appendix, and positive in the other part of the interior of the nonnegativity region.   By using the Sturm sequences, we present a necessary and sufficient condition for a quintic polynomial to be positive and nonnegative  for all positive reals.   We show that for polynomials of a fixed even degree higher than or equal to four, if they have no real roots, then their discriminants take the same sign, which depends upon that degree only, except on an appendix set of dimension lower by two, where the discriminants attain zero.


%15A69:LINEAR AND MULTILINEAR ALGEBRA; MATRIX THEORY- Multilinear algebra, tensor products
%53A45:Classical differential geometry-Vector and tensor analysis
%47A05: Operator Theory-General (adjoints, conjugates, products, inverses, domains, ranges,
%etc.)
%53C35:Global differential geometry- Symmetric spaces
\end{abstract}

%\newpage

\section{Introduction}

In 1988, Schmidt and He\ss \cite{SH88} presented a necessary and sufficient condition for a cubic polynomial to be nonnegative for all positive reals.   In 1994, Ulrich and Watson \cite{UW94} presented a necessary and sufficient condition for a quartic polynomial to be nonnegative for all positive reals.
On the other hand, necessary and sufficient conditions for a quartic polynomial to be positive for all reals have been studied by Gadem and Li \cite{GL64}, Ku \cite{Ku65}, Jury and Mansour \cite{JM81}, Wang and Qi \cite{WQ05}, and Gao \cite{Ga20}.

In this paper, we first present a necessary and sufficient condition for a cubic polynomial to be positive for all positive reals.  We identify the set where the cubic polynomial is nonnegative but not all positive for all positive reals, and explicitly give the points where the cubic polynomial attains zero.

We then reformulate the result of Ulrich and Watson \cite{UW94}.  Based upon this, we derive a necessary and sufficient condition for a quartic polynomial to be nonnegative and positive for all reals.  Comparing with the existing results in the literature, our condition has three merits.   First, it explicitly states a necessary condition.  Second, it explicitly states a symmetric relation between two parameters.  Third, its geometrical meaning is clear.   We also present a theorem on the geometrical properties of the nonnegativity region and the positivity region.   In the interior of the nonnegativity region for all reals, there is an appendix curve.   The discriminant is zero at the appendix, and positive in the other part of the interior of the nonnegativity region.

By using the Sturm sequences, we present a necessary and sufficient condition for a quintic polynomial to be positive and nonnegative  for all positive reals.   This is the first for quintic polynomials.   It also indicates that such conditions are possible for polynomials of even higher degrees.

Finally, we show that such an appendix exists for all even degree polynomials with their degrees higher than or equal to $4$.

These results are useful in automatic control \cite{GL64, Ku65, JM81} and determining copositivity and strict copositivity of symmetric tensors, and positive semidefiniteness of even order symmetric tensors \cite{LZHQ19, LS19, NYZ18, SQ15, SQ20}.  These properties are further useful in optimization \cite{QCC18, QL17}, hypergraphs \cite{CHQ17, CHQ18, Qi13} and  physics \cite{FI19, IKM18, Ka16, Ka18}.

Some preliminary knowledge on the Sturm theorem and quadratic polynomials are given in the next section.
Then, four sections are about cubic polynomials, quartic polynomials, quintic polynomials and higher even degree polynomials, respectively.

All the polynomials studied in this paper are of real coefficients.

\section{Preliminaries}

\subsection{The Sturm Theorem}

Suppose that we have an $m$th degree polynomial $g(t)$.
What is the analytically expressed necessary and sufficient condition such that
$g(t) > 0$ for all $t \ge 0$ ($t > 0$)?

We use the Sturm theorem in classical algebra \cite[pp.52-57]{BPM06} to answer this question.

We may construct the Sturm sequence $\{ g_0, g_1, \cdots, g_l \}$.   Let $g_0 = g$, $g_1 = g'$, and $g_k = -$ rem$(g_{k-1}, g_k)$ for $k \ge 1$, where $g'$ is the derivative of $g$, rem$(g_{k-1}, g_k)$ is the remainder of the division of $g_{k-1}$ by $g_k$.  We have $g_l \not = 0$ and $g_{l+1} = 0$. Then the length $l+1$ of the Sturm sequence is not greater than $m+1$.   The number of variations of the Sturm sequence $S=\{ g_0, g_1, \cdots, g_l \}$ at a real number $\xi$ is the number of sign changes (ignore zero) of the real number sequence
$\{ g_0(\xi), g_1(\xi), \cdots \}$.   Denote it as $V(\xi)$.   By the Sturm theorem, the number of real distinct roots of $g(t) = 0$, for $t \ge 0$, is $V(0) - V(\infty)$.   If $V(0) - V(\infty) = 0$, then $g$ has no positive roots, i.e., $g(t) > 0$ for all $t \ge 0$.    This answers the question.

Note that $g_l$ is the greatest common divisor (GCD) of $g$ and $g'$.  If $g_l$ is not a constant term, then $g$ has multiple roots.   Any root of $g_l$ with multiple $m$ is a $m+1$ multiple root of $g$.

In the following, we simply say that $g$ is strictly copositive if $g(t) > 0$ for all $t \ge 0$.

\subsection{Quadratic Polynomials}

The following result should be known several centuries ago, and is easy to be derived.

\begin{theorem} \label{tt00}
Suppose that we have a quadratic polynomial
$$\phi(t) = t^2 + ut + v,$$
where $v \not = 0$.   Then $\phi(t) > 0$ for all $t \ge 0$ if and only if either (i) $u \ge 0$ and $v > 0$; or (ii) $u < 0$ and $4v > u^2$.    If $u < 0$ and $4v = u^2$, then $\phi(t) \ge 0$ for all $t \ge 0$, and this is the only case that $\phi(t) \ge 0$ for all $t \ge 0$, but not $\phi(t) > 0$ for all $t \ge 0$.

Furthermore, $\phi(t) > 0$ for all $t$ if and only if $4v > u^2$.   If $4v = u^2$, then $\phi(t) \ge 0$ for all $t$, and this is the only case that $\phi(t) \ge 0$ for all $t$, but not $\phi(t) > 0$ for all $t$.
\end{theorem}

\section{Cubic Polynomials}

\subsection{Nonnegaivity Conditions for Non-Degenerate Cubic Polynomials}

 We call a cubic polynomial non-degenerate if its constant is nonzero.  For a non-degenerate cubic polynomial, we may write it as
 \begin{equation} \label{e00}
h(t) = t^3 + pt^2 + qt + r,
\end{equation}
where $r \not = 0$.   If $r < 0$, then it cannot be always nonnegative for all positive reals.  Thus, we may assume that $r > 0$.

By \cite{SH88}, we have the following theorem.

\begin{theorem}  \label{tt01}
Suppose that $h(t)$ is defined by (\ref{e00}).  Then
$h(t) \ge 0$ for all $t \ge 0$ in and only in the following two cases: (A) $p \ge 0$ and $q \ge 0$; (B) (D) $\Delta(h) \le 0$.
Here,
\begin{equation} \label{e02}
\Delta(h) = p^2q^2 + 18pqr - 27r^2 -4p^3r - 4q^3
\end{equation}
is the discriminant of $h(t)$.
\end{theorem}

\subsection{Positivity Conditions for Non-Degenerate Cubic Polynomials}

Suppose that we have a cubic polynomial $h(t)$ defined by (\ref{e00}), with $r> 0$.

Then
$$h(t) = \phi(t)t + r,$$
where $\phi(t) = t^2 + pt + q$.
By Theorem \ref{tt00}, if either (I) $p \ge 0$ and $q \ge 0$, or (II) $q \ge {p^2 \over 4}$, then
$\phi(t) \ge 0$ for all $t \ge 0$.   This implies that $h(t) > 0$ for all $t \ge 0$, as $r > 0$.

We now study the case which is not included in (I) and (II).

Then $q < {p^2 \over 4}$.  We now construct the Sturm sequence for $h(t)$.   Let
$$h_0(t) = h(t) = t^3 + pt^2 + qt + r,$$
$$h_1(t) = h'(t) = 3t^2 + 2pt + q.$$
Then
$$h_0(t) - {t \over 3}h_1(t) = {p \over 3}t^2 + {2q \over 3}t + r,$$
$$h_0(t) - {t \over 3}h_1(t) - {p \over 9}h_1(t) = \left({2q \over 3} - {2p^2 \over 9}\right)t + r - {pq \over 9}.$$
We have
$$h_2(t) = q_2t+r_2,$$
where
\begin{equation} \label{e01}
q_2 = {2p^2 \over 9} - {2q \over 3}, \ r_2 = {pq \over 9} -r.
\end{equation}

Since $q < {p^2 \over 4}$, we have $q_2 > 0$.

Then
$$h_3(t) \equiv r_3 = - h_1\left(-{r_2 \over q_2}\right) = -{3r_2^2 \over q_2^2} + {2pr_2 \over q_2} -q = {\Delta(h) \over 9q_2^2},$$
where
$$\Delta(h) \equiv p^2q^2 + 18pqr - 27r^2 -4p^3r - 4q^3 = -81r_2^2 + 54pq_2r_2 -27qq_2^2.$$
The expression of $\Delta(h)$ can be found in \cite{GKZ94}.

Suppose that $\Delta(h) = 0$.   Then $h_2(t)$ is the GCD of $g$ and $g'$,
$$h(t) = (t-\alpha)^2(t+\beta),$$
where
$$\alpha = -{r_2 \over q_2}.$$
We have $q_2 > 0$.   We now show that $r_2 < 0$.  Assume that $r_2 \ge 0$.   This means that $pq \ge 9r$.
Since we have already exclude the case that $p \ge 0$ and $q \ge 0$, this means $p < 0$, $q < 0$ and $pq \le 9r$.   These further implies that $\Delta(h) < 0$, contradicting the assumption that $\Delta(h) = 0$.

We have $\alpha^2\beta = r > 0$.   Thus, $\beta > 0$ and $\alpha \not = 0$.   Therefore, $h(t) \ge 0$ for all $t \ge 0$.   Furthermore, in this case, we always have $t = \alpha > 0$ such that $h(t) = 0$.

Suppose now $\Delta(h) \not = 0$.   We have
$$S = \{ h_0, h_1, h_2, h_3 \},$$
$$S(\infty) = \{ 1, 1, q_2, \Delta(h) \},$$
$$S(0) = \{ 1, q, r_2, \Delta(h) \}.$$

Assume that $\Delta(h) > 0$.   Then $V(\infty) = 0$.
Thus, $h(t) > 0$ for all $t \ge 0$ if and only if $V(0) = 0$, i.e., $q \ge 0$ and $r_2 \ge 0$.   However, this implies that $p > 0$ and $q > 0$, which has already been covered by Case (I).

Assume that $\Delta(h) < 0$.   Then $V(\infty) = 1$.
Thus, $h(t) > 0$ for all $t \ge 0$ if and only if either (i) $r_2 \le 0$, or (ii) $q \ge 0$ and $r_2 \ge 0$.   However, Case (ii) has already been covered by Case (I) as discussed above.  Thus, only Case (i) needs to be considered.    As discussed above,  the case that $r_2 > 0$, i.e., $pq > 9r$, is covered by the condition (I) $p \ge 0$ and $q \ge 0$, and (III) $\Delta(h) < 0$.    Thus, the condition that $r_2 \le 0$ is not necessary here.

Thus, we have one more case such that $h(t) > 0$ for all $t \ge 0$: (III) $\Delta(h) < 0$.

Finally, Case (II) $q \ge {p^2 \over 4}$ is coved by the union of Cases (I) and (III).

Thus, we have the following theorem.

\begin{theorem}  \label{tt1}
Suppose that $h(t)$ is defined by (\ref{e00}).  Then
$h(t) > 0$ for all $t \ge 0$ in and only in the following two cases: (A) $p \ge 0$ and $q \ge 0$; (B) $\Delta(h) < 0$.
Here, $\Delta(h)$ is the discriminant, defined by (\ref{e02}).

If $\Delta(h) = 0$ and either $p < 0$ or $q < 0$, then $h(t) \ge 0$ for all $t \ge 0$, but there is
$$\alpha = {9r-pq \over 2p^2 - 6q} > 0$$
such that $h(\alpha) = 0$.

In all the other cases, there is a $t> 0$ such that $h(t) < 0$.
\end{theorem}

%{\bf Remark} For quadratic and cubic polynomials, if the discriminant is negative, then the polynomial is positive for all positive reals.   This is no longer true for quartic polynomials.   Consider the following quartic polynomial
%$$f(t) = t^4 - 2t^3 +t^2 - 2t + 1.$$
%It is negative for $t=1$.  But by the formula in Theorem 2 of \cite{UW94}, we have $\Delta(f) = -357 < 0$.
%Hence, the positivity condition for polynomials of higher degrees would be much completed.

\section{Quartic Polynomials}

\subsection{Nonnegativity Conditions}

Suppose that we have a quartic polynomial
\begin{equation} \label{qe1}
f(t) = t^4 + \alpha t^3 + \beta t^2 + \gamma t + 1.
\end{equation}
By \cite{UW94}, its discriminant has the form:
\begin{equation} \label{qe2}
\Delta(f) = 4[\beta^2-3\alpha \gamma +12]^3-[72\beta +9\alpha \beta \gamma-2\beta^3-27\alpha^2-27\gamma^2]^2.
\end{equation}

By \cite{UW94}, we have the following theorem.

\begin{theorem} \label{qt1}
Suppose that $f$ is defined by (\ref{qe1}).   Let
%\begin{equation} \label
$$\Lambda_1 = (\alpha -\gamma )^2-16(\alpha +\beta +\gamma +2),\ \Lambda_2 = (\alpha -\gamma )^2 -{4(\beta +2) \over \sqrt{\beta -2}}(\alpha +\gamma  + 4\sqrt{\beta -2}).$$
%\end{equation}
Then $f(t) \ge 0$ for all $t > 0$ if and only if either

(1) $\beta < -2$, $\Delta(f) \le 0$ and $\alpha +\gamma > 0$; or

(2) $-2 \le \beta \le 6$ and either (i) $\Delta(f) \le 0$ and $\alpha + \gamma  > 0$, or (ii) $\Delta(f) \ge 0$ and $\Lambda_1 \le 0$; or

(3) $\beta > 6$ and either (i) $\Delta(f) \le 0$ and $\alpha +\gamma > 0$, or (ii) $\alpha > 0$ and $\gamma > 0$; or (iii) $\Delta(f) \ge 0$ and $\Lambda_2 \le 0$.
\end{theorem}

We may combine (1), (2i) and (3i) as

(A) $\Delta(f) \le 0$ and $\alpha +\gamma  > 0$.

On the other hand, as
$$f(t) = (t^4 + \beta t^2 + 1) + \alpha t^3 + \gamma t,$$
if $\alpha \ge 0$, $\gamma \ge 0$ and $\beta \ge -2$, then $f(t) \ge 0$ for all $t \ge 0$.    Thus, we have

(B) $\alpha \ge 0$, $\gamma \ge 0$ and $\beta \ge -2$.

Then (B) covers (3ii).

We may replace the constraints $\Lambda_1 \le 0$ and $\Lambda_2 \le 0$, by some constraints which are linear with respect to $\alpha$ and $\gamma$.

For the condition $\Lambda_1 \le 0$, from the discussion in Sections 2 and 3 of \cite{UW94}, we may replace it by $|\alpha - \gamma | \le r,$
where $r$ is the value of $|\alpha - \gamma |$ at the two intersection points of $\alpha + \gamma = 0$ and
the curve $\Gamma_\beta$ in \cite{UW94}.   These are two points where $\alpha(t) + \gamma(t) = 0$, for $\alpha(t)$ and $\gamma(t)$ given by (7a) and (7b) of \cite{UW94}.   Solving this, we find that $r = 4\sqrt{\beta +2}$.    Thus, we may use $|\alpha - \gamma| \le 4\sqrt{\beta+2}$ to replace $\Lambda_1 \le 0$.   Note that we may rewrite the constraint $|\alpha - \gamma| \le 4\sqrt{\beta+2}$ such that
$\alpha$ and $\gamma$ are linear there.

For the condition $\Lambda_2 \le 0$, by Figure 1 of \cite{UW94}, we may replace it by $|\alpha - \gamma| \le 4\sqrt{\beta+2}$ and
$\alpha + \gamma \ge -4\sqrt{\beta-2}$.

We have the following theorem.

\begin{theorem} \label{qt2}
Suppose that $g(z) = az^4 + bz^3 + cz^2 + dz + e$ be a quartic polynomial with real coefficients and $a > 0$ and $e > 0$.   Let
$$\alpha = ba^{-3 \over 4}e^{-1 \over 4},\ \beta = ca^{-1 \over 2}e^{-1 \over 2},\ \gamma = da^{-1 \over 4}e^{-3 \over 4},$$
$f$ is defined by (\ref{qe1}), and $\Delta(f)$ is defined by (\ref{qe2}).
Then $f(t) \ge 0$ for all $t > 0$ if and only if either

(A) $\Delta(f) \le 0$ and $\alpha + \gamma > 0$; or

(B) $\alpha \ge 0$, $\gamma \ge 0$ and $\beta \ge -2$; or

(C) $\Delta(f) \ge 0$, $|\alpha - \gamma | \le 4\sqrt{\beta+2}$ and either (i) $-2 \le \beta \le 6$, or (ii) $\beta > 6$ and $\alpha + \gamma \ge -4\sqrt{\beta-2}$.

%Furthermore, $f(t) > 0$ for all $t > 0$ if and only if either

%(A') $\Delta(f) < 0$ and $\alpha+\gamma > 0$; or

%(B') $\alpha \ge 0$, $\gamma \ge 0$ and $\beta \ge -2$ except the point that $\alpha = 0$, $\gamma = 0$ and $\beta = -2$; or

%(C') $\Delta(f) > 0$, $|\alpha - \gamma| \le 2\sqrt{\beta+2}$ and either (i) $-2 \le \beta \le 6$, or (ii) $\beta > 6$ and  $\alpha + \gamma \ge -4\sqrt{\beta-2}$.

\end{theorem}
{\bf Proof} We may convert the general quartic polynomial $g(z)$ to $f(t)$, defined by (\ref{qe1}), as in \cite{UW94}.

We have already discussed to convert Conditions (1), (2) and (3) of Theorem \ref{qt1} to Conditions (A), (B) and (C).  %Conditions (A'), (B') and (C') describe the interior of the positivity region.  In the interior of the positivity region, $\Delta(f) \not = 0$.   On the other hand, if $f(t) \ge 0$ for all $t \ge 0$, but there is a point $t^* > 0$ such that $f(t^*) = 0$, then $t^*$ must be a multiple root of $f$.   We have $\Delta(f) = 0$.   Thus, under conditions (A'), (B') and (C'), $f(t) > 0$ for all $t \ge 0$.   If $f$ satisfies conditions (A), (B) and (C), but dies not satisfy conditions (A'), (B') and (C'), then there are quartic polynomials $f_{\epsilon}$ such that for each $\epsilon$, there is a $t_\epsilon > 0$ such that $f_\epsilon(t_\epsilon) = 0$, and $f_\epsilon \to f$ as $\epsilon \to 0$.  By boundedness of the roots of $f(t)$, $t_\epsilon$ is bounded for $\epsilon \to 0$. Then $t_\epsilon$ has a limiting point $t^* \ge 0$.   We have $f(t^*) = 0$.
\qed

A merit of the format of Theorem \ref{qt2} is that it is somewhat convenient to derive nonnegativity conditions for all $t$ from conditions (A), (B) and (C).

\subsection{Nonnegativity and Positivity Conditions for All Reals}

We now consider the conditions that $f(t) \ge 0$ ($f(t) > 0$) for all $t$.  This is equivalent to the condition that $f(t) \ge 0$ ($f(t) > 0$) and $g(t) \ge 0$ ($g(t) > 0$) for all $t$, where $g(t)$ is defined by
$$g(t) = t^4 - \alpha t^3 + \beta t^2 - \gamma t + 1.$$
Note that $\Delta(f) = \Delta(g)$.  By Theorem \ref{qt2}, we have the following theorem.

\begin{theorem} \label{qt3}
Suppose that $g(z) = az^4 + bz^3 + cz^2 + dz + e$ is a quartic polynomial with real coefficients and $a > 0$ and $e > 0$.   Let
$$\alpha = ba^{-3 \over 4}e^{-1 \over 4},\ \beta = ca^{-1 \over 2}e^{-1 \over 2},\ \gamma = da^{-1 \over 4}e^{-3 \over 4},$$
$f$ be defined by (\ref{qe1}), and $\Delta(f)$ be defined by (\ref{qe2}).
Then $f(t) \ge 0$ for all $t$ if and only if $\Delta(f) \ge 0$, $|\alpha - \gamma| \le 4\sqrt{\beta+2}$ and either (i) $-2 \le \beta \le 6$, or (ii) $\beta > 6$ and  $|\alpha + \gamma| \le 4\sqrt{\beta-2}$.

Furthermore, $f(t) > 0$ for all $t$ if and only if either

(A) $\Delta(f) = 0$, $\alpha = \gamma$, $\alpha^2 +8 = 4\beta < 24$; or

(B) $\Delta(f) > 0$, $|\alpha - \gamma| \le 4\sqrt{\beta+2}$ and either (i) $-2 \le \beta \le 6$, or (ii) $\beta > 6$ and $|\alpha + \gamma| \le 4\sqrt{\beta-2}$.
\end{theorem}
{\bf Proof}   Let $\bar f(t) = f(-t)$.  Then
$$\bar f(t) = t^4 -\alpha t^3 + \beta t^2 - \gamma t + 1.$$
Then $f(t) \ge 0$ for all $t$ if and only if $f(t) \ge 0$ and $\bar f(t) \ge 0$ for all $t \ge 0$.
From Theorem \ref{qt2}, we have conditions that $f(t) \ge 0$ and $\bar f(t) \ge 0$ for all $t \ge 0$.
Combining these conditions, i.e., taking the intersection of the conditions that $f(t) \ge 0$ for all $t \ge 0$, and the conditions that $\bar f(t) \ge 0$ for all $t \ge 0$, we have the conditions for the first part of this theorem.

We now prove the second part of this theorem.  If $\Delta(f) > 0$, it is in the interior of the nonnegativity region, we have $f(t) > 0$ for all $t$.  We have condition (B).  Now assume that $\Delta(f) =0$.   Then $f(t)$ has a multiple root.  To make $f(t) > 0$ for all $t$, this multiple toot has to be complex.   Thus, $f(t) = (t^2 + ut + 1)^2$ with $u^2 < 4$.  Comparing the coefficients of
$$f(t) = (t^2+ ut+ 1)^2 = t^4 + \alpha t^3 + \beta t^2 + \gamma t +1,$$
We have $\alpha = \gamma = 2u$ and $u^2 + 2 = \beta$.   Then,
we have condition (A).
%Like the proof of Theorem \ref{qt2}, considering the interior of the positivity region, we have the second part of this theorem.
\qed

Geometrically, for fixed $\beta$, the nonnegativity region of $(\alpha, \gamma)$ for $g(t)$ is the convex hull of the central convex compact part of the positive sign region in Figure 5b and 5c of \cite{UW94}, while the positivity region is the interior of the nonnegativity region.   In Figures 5b of \cite{UW94},
for $2 \le \beta < 6$, there are two points missing.  These two points $(\alpha, \beta, \gamma)$ are defined by $\alpha = \gamma$ and $\alpha^2 + 8 = 4\beta$.   At these two points, $\Delta(f) = 0$.  These two points are missing in Figure 5b of \cite{UW94}, but do not affect the final results of \cite{UW94}.

% and the intersection of the boundary with the plane $\alpha = \gamma$ and $2 \le \beta < 6$.

Note that $\beta +2 \ge 0$ is a necessary condition such that $f(t) \ge 0$ for all $t$.   This explicitly stated  in Theorem \ref{qt3}.  Actually, if $f(t) \ge 0$ for all $t$, then
$$\beta + 2 = {f(1) + f(-1) \over 2} \ge 0.$$

The following theorem presents the geometrical features of the nonnegativity region and the positivity region.

\begin{theorem} \label{qgt}
Suppose that $f(t)$ is defined by (\ref{qe1}).  Let
$$S = \{ (\alpha, \beta, \gamma) : f(t) > 0 \ {\rm for\ all}\ t \}$$
and
$$\bar S = \{ (\alpha, \beta, \gamma) : f(t) \ge 0 \ {\rm for\ all}\ t \}.$$
Then $\bar S$ is a closed convex cone with a recession direction $(0, 1, 0)$ and an apex $(0, -2, 0)$.  It is symmetric with respect to $\alpha$ and $\gamma$, i.e., if $(\alpha, \beta, \gamma) \in \bar S$, then
$(\gamma, \beta, \alpha), (-\alpha, \beta, -\gamma) \in \bar S$.
For any $(\alpha, \beta, \gamma)$ in the interior of $\bar S$, if $2 \le \beta < 6$, $\alpha = \beta$ and $\alpha^2 + 8 = 4\beta$, we have $\Delta(f) = 0$.  Otherwise, we have $\Delta(f) > 0$.   For any
$(\alpha, \beta, \gamma)$ at the boundary of $\bar S$, we have $\Delta(f) = 0$.  The set $S$ is
the interior of $\bar S$. %and the intersection of the boundary of $\bar S$ with the plane $\alpha = \gamma$ and $2 \le \beta < 6$.
\end{theorem}
{\bf Proof}  Consider $\bar S$.  Suppose that $f(t)$ and $g(t)$ are two quartic polynomials with the form of (\ref{qe1}),  $f(t) \ge 0$ and $g(t) \ge 0$ for all $t$.  Then ${1 \over 2}(f(t)+g(t))$ is still a quartic polynomial with the form of (\ref{qe1}), and is nonnegative for all $t$.   This shows that $\bar S$ is convex.  Similarly, we may show that it is closed.   Let $f(t)$ be defined by (\ref{qe1}), and $\delta > 0$.  Let
$$g(t) = t^4 + \alpha t^3 +(\beta + \delta)t^2 + \gamma t + 1.$$
Then $g(t) = f(t) + \delta t^2 \ge f(t) \ge 0$ for all $t$.   This shows that $(0, 1, 0)$ is a recession direction of $\bar S$ and $\bar S$ is a cone.  From Theorem \ref{qt3}, the point $(\alpha, \beta, \gamma) = (0, -2, 0)$ is the unique point of $\bar S$ with the smallest value of $\beta$.   Hence, it is an apex of $\bar S$.  The other properties of $\bar S$ and $S$ can be derived from Theorem \ref{qt3} accordingly.
\qed

\includegraphics[width=10.5cm]{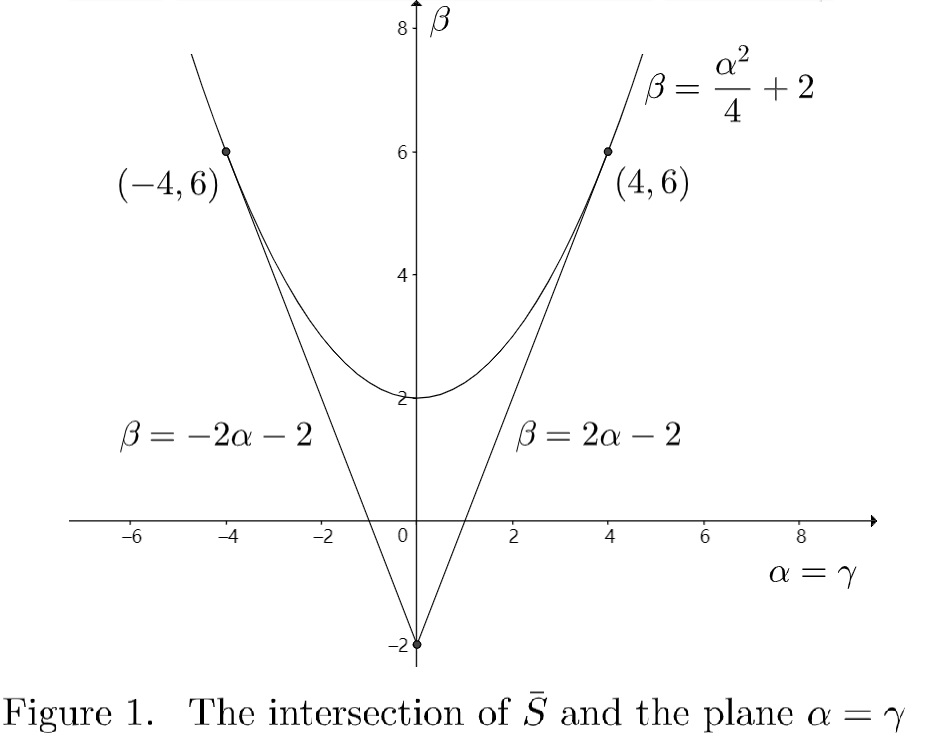}

In Figure 1, the $\alpha = \gamma$ plane is exhibited.   The set $\bar S$ is the area above $\beta = 2|\alpha|-2$ for $|\alpha| \le 4$, and the area above the parabola $\beta = {\alpha^2 \over 4} + 2$ for $|\alpha| \ge 4$.  The discriminant $\Delta(f)$ vanishes at the line segments $\beta = 2|\alpha| - 2$ for $|\alpha| \le 4$ and the parabola $\beta = {\alpha^2 \over 4} + 2$ for all $\alpha$.   The truncated parabola $\beta = {\alpha^2 \over 4} + 2$ for $|\alpha| < 4$ is in the interior of $\bar S$.    In Figure 2, the $\alpha = -\gamma$ plane is exhibited.  The set $\bar S$ is the area above the parabola $\beta = {\alpha^2 \over 4} - 2$.   We may see that $(\alpha, \beta, \gamma) = (0, -2, 0)$ is the apex of the cone $\bar S$.

The truncated parabola $\beta = {\alpha^2 \over 4} + 2$ for $|\alpha| < 4$ in Figure 1 is a special part of
the surface $\Delta(f) = 0$.   The other part of the surface $\Delta(f) = 0$ is of dimension $2$, but the
truncated parabola  $\beta = {\alpha^2 \over 4} + 2$ for $|\alpha| < 4$ is of dimension $1$.   Thus, we call it the {\bf appendix} of the surface $\Delta(f) = 0$.  Such an appendix exists for polynomials with their degrees higher than or equal to $4$.

\includegraphics[width=10.5cm]{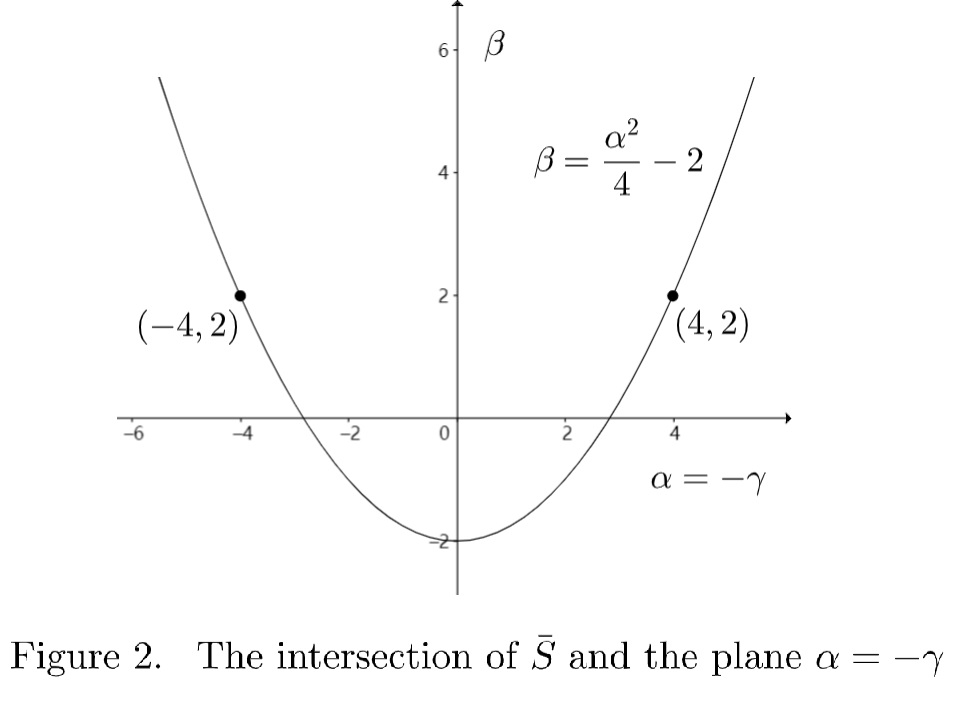}

This problem has been studied by Gadem and Li \cite{GL64}, Ku \cite{Ku65}, Jury and Mansour \cite{JM81}, Wang and Qi \cite{WQ05}, and Gao \cite{Ga20}.

We may compare Theorem \ref{qt3} with the result of \cite{WQ05}, which is a correction of the result of \cite{Ku65}.    The following is the result of \cite{WQ05}.

The polynomial considered in \cite{WQ05} has the following form:
$$f(t) = a_0t^4 + 4a_1t^3 + 6a_2t^2 + 4a_3t + a_4.$$
Thus,
$$a_0 = a_4 = 1,\ a_1 = {\alpha \over 4},\ a_2 = {\beta \over 6},\ a_3 = {\gamma \over 4}.$$
Let
$$G = a_0^2a_3 -3a_0a_1a_2 + 2a_1^3 = {\gamma \over 4} - {\alpha\beta \over 8} + {\alpha^3 \over 32}= {8\gamma -4\alpha\beta + \alpha^3 \over 32}.$$
$$H = a_0a_2 - a_1^2 = {\beta \over 6} - {\alpha^2 \over 16} = {8\beta - 3\alpha^2 \over 48}.$$
$$I = a_0a_4-4a_1a_3 + 3a_2^2 = 1-{\alpha\gamma \over 4} + {\beta^2 \over 12}= {12 - 3\alpha\gamma + \beta^2 \over 12}.$$
$$J = a_0a_2a_4 + 2a_1a_2a_3-a_1^2a_4-a_0a_3^2 - a_2^3 = {\beta \over 6} + {\alpha\beta\gamma \over 48} -{\alpha^2+\gamma^2 \over 16} - {\beta^3 \over 216} = {72\beta + 9\alpha\beta\gamma - 27\alpha^2 - 27\gamma^2 -2\beta^3 \over 432}.$$
$$\Delta = I^3 -27J^3 = \left({12 - 3\alpha\gamma + \beta^2 \over 12}\right)^3 -{(72\beta + 9\alpha\beta\gamma - 27\alpha^2 - 27\gamma^2 -2\beta^3)^2 \over 3 \times 36^2} = {\Delta(f) \over 4 \times 12^3}.$$
Hence, $\Delta$ in \cite{WQ05} has the same sign as $\Delta(f)$ and plays the same role as $\Delta(f)$.
The result of \cite{WQ05} can be formulated with the language of this paper as follows.

\begin{theorem} \label{qt4}
Let $f(t)$ be defined by (\ref{qe1}).  Then $f(t) > 0$ for all $t$ if and only if either

(1) $\Delta(f) = 0$, $G = 0$, $12H^2-I = 0$ and $H > 0$; or

(2) $\Delta(f) > 0$ and either (i) $H \ge 0$, or (ii) $H < 0$ and $12H^2-I < 0$.
\end{theorem}

We have
$$12H^2 - I = 12\left({8\beta - 3\alpha^2 \over 48}\right)^2 - {12 - 3\alpha\gamma + \beta^2 \over 12} = {64\beta^2 -48\alpha^2\beta+9\alpha^4 - 192 +48\alpha\gamma -16\beta^2 \over 192}$$
$$= {16\beta^2 -16\alpha^2\beta+3\alpha^4 - 64 +16\alpha\gamma \over 64}.$$

\begin{corollary} \label{qc1}
Let $f(t)$ be defined by (\ref{qe1}).  Then $f(t) \ge 0$ for all $t$ if and only if
$\Delta(f) \ge 0$ and either (i) $H \ge 0$, or (ii) $H < 0$ and $12H^2-I \le 0$.
\end{corollary}

The conclusions of Theorem \ref{qt3} should be the same with the conditions of Theorem \ref{qt4} and Corollary \ref{qc1}.   However, Theorem \ref{qt3} has three merits.   First, it explicitly stated that a necessary condition is that $\beta \ge -2$.  Second, it treats $\alpha$ and $\gamma$ in a symmetric way.
Third, its geometrical meaning is clear, as shown by Theorem \ref{qgt}.  The set indicated by Theorem \ref{qt3} (A) and Theorem \ref{qt4} (1) is the appendix of the surface $\Delta(f) = 0$.  It is actually in the interior of the nonnegativity region.
%We should check this.
%First, we should see that $\beta \ge -2$ is implicitly implied by the conditions of Corollary \ref{qc1}.
%Actually, $H \ge 0$ implies that $\beta \ge 0$.  Now we wish to check if $12H^2-I \le 0$ implies $\beta \ge -2$.
%$$16\beta^2 -16\alpha^2\beta+3\alpha^4 - 64 +16\alpha\gamma = 16\left(\beta - {\alpha^2 \over 4}\right)^2 - \alpha^4 - 64 +16\alpha\gamma.$$

\section{Quintic Polynoimals}

\subsection{Quintic Polynoimals with Multiple Roots}

Suppose that we have a quintic polynomial
\begin{equation} \label{e1}
g(t) = t^5 + at^4 + bt^3 + ct^2 + dt + e,
\end{equation}
where $e > 0$.

As we discussed in Section 2, by construct the Sturm sequence of $g$, we may find the GCD of $g$ and $g'$ if $g$ has a multiple root.

\medskip

(A). Suppose that the GCD of $g$ and $g'$ is a linear polynomial. We may assume that it is $\phi(t) = t - \alpha$.   Then $g$ has a double root $\alpha$.  Since $e > 0$, $\alpha \not = 0$.  We may write that
\begin{equation} \label{e03}
g(t) = h(t)(t-\alpha)^2,
\end{equation}
where $h(t) = t^3 + pt^2 + qt + r$.  Then
\begin{equation} \label{e030}
p = a + 2\alpha,\ q = b+2\alpha p + \alpha^2,\ r = {e \over \alpha^2}.
\end{equation}
Since $e > 0$, $r > 0$, i.e., $h$ has the form (\ref{e00}).  By Theorem \ref{tt1}, we have the following proposition.

\begin{proposition} \label{p1}
Suppose that $g(t)$ is defined by (\ref{e1}), $g$ and $g'$ has a GCD $\phi(t) = t - \alpha$.  Then $g$ has the form (\ref{e03}), where $h(t)$ has the form (\ref{e00}), with $p, q$ and $r$, given by (\ref{e030}).   Then $g(t) \ge 0$ for all $t \ge 0$, if and only if $h(t) \ge 0$ for all $t \ge 0$.   If furthermore $\alpha < 0$, then $g(t) > 0$ for all $t \ge 0$, if and only if $h(t) > 0$ for all $t \ge 0$.
\end{proposition}

\medskip

(B). Suppose that the GCD of $g$ and $g'$ is a quadratic polynomial. We may assume that it is $\phi(t) = t^2 + ut + v$.   There are three subcases.

(B1). $u^2 \not = 4v$.   In this subcase, $\phi(t)$ has two distinct roots $\alpha$ and $\beta$.  Then $g$ has a double root $\alpha$ and a double root $\beta$.  Since $e > 0$, $\alpha \not = 0$ and $\beta \not = 0$.   We have
\begin{equation} \label{e04}
g(t) = \phi(t)^2(t-\gamma) = (t-\alpha)^2(t-\beta)^2(t-\gamma),
\end{equation}
where $\gamma$ is the fifth root of $g$.   Since $e > 0$, $\gamma < 0$.   Then we have the following conclusion.

\begin{proposition} \label{p2}
Suppose that $g(t)$ is defined by (\ref{e1}), $g$ and $g'$ has a GCD $\phi(t) = t^2 + ut + v$.  If $u^2 \not = 4v$, then $g(t) \ge 0$ for all $t \ge 0$.   If furthermore $\phi(t) > 0$ for all $t \ge 0$, then $g(t) > 0$ for all $t \ge 0$.
\end{proposition}

We may use Theorem \ref{tt00} to determine if $\phi(t) > 0$ for all $t > 0$.

(B2). $u^2 = 4v$. Then $\phi$ has a double root $\alpha = -{u \over 2}$, and $\alpha$ is a triple root of $g$.   Since $e > 0$, $\alpha \not = 0$.  We have
\begin{equation} \label{e04}
g(t) = \left(t+ {u \over 2}\right)^3\psi(t),
\end{equation}
where
\begin{equation} \label{e05}
\psi(t) = t^2 + \hat ut + \hat v.
\end{equation}
Then
\begin{equation} \label{e06}
\hat v = {8e \over u^3}, \ \hat u = a - {3u \over 2}.
\end{equation}
We have the following conclusion.

\begin{proposition} \label{p3}
Suppose that $g(t)$ is defined by (\ref{e1}), $g$ and $g'$ has a GCD $\phi(t) = t^2 + ut + v$.  Suppose that $u^2 = 4v$.  Then $u \not = 0$.  Let $\psi(t)$ be calculated by (\ref{e05}) and (\ref{e06}).  Then $g(t) \ge 0$ for all $t \ge 0$ if and only if $u > 0$ and $\psi(t) \ge 0$ for all $t \ge 0$,  and $g(t) > 0$ for all $t \ge 0$ if and only if $u > 0$ and $\psi(t) > 0$ for all $t \ge 0$.
\end{proposition}

We may use Theorem \ref{tt00} to determine the situation of $\psi(t)$.

(C). Suppose that the GCD of $g$ and $g'$ is a cubic polynomial. We may assume that it is $h(t) = t^3 + pt^2 + qt + r = (t-\alpha)^2(t+\beta)$.   There are two subcases.

(C1) $q = {p^2 \over 3}$ and $r = {p^3 \over 27}$.  This implies $\alpha = -\beta = -{p \over 3}$, and
$$g(t) = \left(t + {p \over 3}\right)^4(t - \gamma).$$
Since $e > 0$, $\gamma < 0$, $p \not = 0$.  Thus, $g(t) \ge 0$ for all $t \ge 0$.   If furthermore $p > 0$, then $g(t) > 0$ for all $t \ge 0$.

(C2) Otherwise, $\alpha \not = -\beta$, $\alpha \not = 0$ and $\beta \not = 0$.   We have $\alpha$ from (\ref{e020}) and (\ref{e01}).  Then
$$\beta = {r \over \alpha^2}.$$
We have
$$g(t) = (t-\alpha)^3(t+\beta)^2.$$
Then $-\alpha^3\beta^2 = e$.   Since $e > 0$, $\alpha < 0$.   Thus, we always have $g(t) \ge 0$ for all $t \ge 0$.   If furthermore $r > 0$, then $g(t) > 0$ for all $t \ge 0$.

We have the following proposition.

\begin{proposition} \label{p4}
Suppose that $g(t)$ is defined by (\ref{e1}), $g$ and $g'$ has a GCD $h(t) = t^3 + pt^2 + qt + r$.
Then $g(t) \ge 0$ for $t \ge 0$.    We have $g(t) > 0$ for all $t \ge 0$ if and only if either
(i) $q = {p^2 \over 3}$, $r = {p^3 \over 27}$ and $p > 0$; or (ii) $r > 0$ and either $q \not = {p^2 \over 3}$ or $r \not = {p^3 \over 27}$.
\end{proposition}

The case that the GCD of $g$ and $g'$ is a quartic polynomial will be analyzed in the next section.  In that case, $g(t) > 0$ for all $t \ge 0$.

\subsection{General Quintic Polynomials}

For constructing the Sturm sequence of a quintic polynomial $g(t)$, we keep to use $a, b, c, d$ and $e$ to denote the coefficients of the polynomials $g_i$.
We have
$$g_0(t) = g(t) = t^5 + at^4 + bt^3 + ct^2 + dt + e,$$
and
$$g_1(t) = g'(t) = 5t^4 + 4at^3 + 3bt^2 + 2ct + d.$$
For $i \ge 2$, we use $b_i, c_i, d_i$ and $e_i$ to denote the coefficients of $g_i(t)$, with $b_i$ for the coefficient of $t^3$, $c_i$ for $t^2$, $d_i$ for $t$ and $e_i$ for the constant term.   If there are different coefficients used, we use an additional index to distinguish them.   For example, in the following, $e_3$ was used in Case (2), and in Case (3), we use $e_{3, 1}$ to denote a new coefficient.  We use $g_3$ to denote different polynomials in different parts of this section, while $e_3$ and $e_{3,1}$ are uniquely used in this section, as $g_3$ will not appear in the theorems in the next sections, but  $e_3$ and $e_{3, 1}$ play a role in establishing those theorems.

However, $e_3$ and $e_{3, 1}$ are fractional functions of the coefficients $a, b, c, d$ and $e$.  This is not convenient.  Then, we use the bar symbol to denote some polynomial functions to replace them.   For example, we write $e_3 = {\bar e_3 \over 25d_2^4}$, where $\bar e_3$ is a polynomial function of the coefficients $a, b, c, d$ and $e$, and we use $\bar e_3$ instead of $e_3$ in our theorems.    Totally, eleven such fractional functions are replaced.

In the following, $S(0)$ and $S(\infty)$ are the Sturm sequence at $0$ and $\infty$.  Since only signs are important, we may replace their entries by other numbers as long as the signs are not changed.

Note that when there are no multiple roots, either $g(t)$ has no positive root, or has at least two positive roots.  This means that either $V(0) - V(\infty) = 0$ or $V(0) - V(\infty) \ge 2$.

%2. Second, in the following, starting from Case (3), after we have the Sturm sequence, we analyze every possible values of $V(0) = V(\infty)$ one by one, such that we will not miss any case.

%2. We use $A \vee B$ to express that conditions $A$ or $B$,  and $A \wedge B$ to express condition $A$ and $B$.

We now construct the Sturm sequence for $g$.
We have
$$g(t) - {t \over 5}g_1(t) = {a \over 5}t^4 + {2b \over 5}t^3 + {3c \over 5}t^2 + {4d \over 5}t + e,$$
$$g(t) - {t \over 5}g_1(t) - {a \over 25}g_1(t) = \left({2b \over 5} - {4a^2 \over 25}\right)t^3 + \left({3c \over 5}- {3ab \over 25}\right)t^2 + \left({4d \over 5}-{2ac \over 25}\right)t + \left(e -{ad \over 25}\right),$$
\begin{eqnarray*}
g_2(t) & = & {4a^2-10b \over 25}t^3 + {3ab-15c \over 25}t^2 + {2ac-20d \over 25}t + {ad-25e \over 25}\\
& = & b_2t^3 + c_2t^2 + d_2t +e_2.
\end{eqnarray*}

There are four possibilities.

(1) $b_2 = c_2 = d_2= 0$.   We further divide this case to two subcases.

(1A) $e_2 = 0$.  This implies that $b = {2a^2 \over 5}$, $c = {ab \over 5}={2a^3 \over 25}$, $d = {ac \over 10} = {a^4 \over 125}$ and $e = {ad \over 25} = {a^5 \over 3125}$.  Thus,
$$g(t) = \left(t+{a \over 5}\right)^5.$$
Since $e > 0$, we have $a > 0$.   Then, $g$ is strictly copositive in this subcase.

(1B) $e_2 \not = 0$.  Then $S = \{ g_0, g_1, g_2 \}$, and $b = {2a^2 \over 5}$, $c = {ab \over 5}$ and $d = {ac \over 10}$, and  $g_2(t) = e_2$.   We have
$$d = {ac \over 10} = {a^2b \over 50} = {a^4 \over 125} \ge 0.$$
Then $V(\infty) = 0$ if $e_2 > 0$.   Otherwise $V(\infty) = 1$.   On the other hand, $V(0) = 1$ if $e_2 < 0$, $V(0) = 0$ if $e_2 > 0$.  Thus, $V(0) - V(\infty) \equiv 0$, and $g$ is strictly copositive in this subcase.

Thus, $g$ is strictly copositive in Case (1).

(2) $b_2 = c_2 = 0$ and $d_2 \not = 0$.  Then $b = {2a^2 \over 5}$, $c = {ab \over 5} = {2a^3 \over 25}$,
$$g_2(t) = d_2t+e_2,$$
$$g_1(t) = 5t^4 + 4at^3 + {6a^2 \over 5}t^2 + {4a^3 \over 25}t + d,$$
$${\rm rem}\{g_1, g_2\}= g_1\left(-{e_2 \over d_2}\right),$$
$$g_3(t) \equiv e_3 = -g_1\left(-{e_2 \over d_2}\right)= {\bar e_3 \over 25d_2^4},$$
where
$$\bar e_3 = -125e_2^4+100ae_2^3d_2-30a^2e_2^2d_2^2+4a^3e_2d_2^3-dd_2^4.$$

We need to divide this case to two subcases.

(2A) $e_3 = 0$, i.e., $\bar e_3 = 0$.  In this subcase, $\alpha = -{e_2 \over d_2}$ is a double root of $g$.  We may apply Proposition \ref{p1} and Theorem \ref{tt1} to determine the situation of $g(t)$.

(2B) $e_3 \not = 0$, i.e., $\bar e_3 \not = 0$.   Then
$$S = \{ g_0, g_1, g_2, g_3 \},$$
$$S(\infty) = \left\{+\infty, +\infty, d_2, e_3 \right\},$$
$$S(0) = \{e, d, e_2, e_3 \}.$$
Clearly, we may replace $e_3$ by $\bar e_3$ here.
Then $V(0) = V(\infty)$ if and only if (i) $V(0) = V(\infty) =0$; or (ii) $V(0) = V(\infty) = 1$; or (iii) $V(\infty) = 2$.   As we discussed early, $V(0) - V(\infty) \ge 2$.   Hence, in (iii), we have $V(0)=2$ too as $V(0)$ cannot be $4$ in this subcase.

Thus, in this subcase, $g$ is strictly copositive if and only if either (i) $\bar e_3 > 0$, $d_2 \ge 0$, $d \ge 0$, $e_2 \ge 0$; or (ii) $\bar e_3 < 0$, and $d \ge 0$ if $e_2 > 0$; or (iii) $\bar e_3 > 0$, $d_2 < 0$.

(3) $b_2 = 0$ and $c_2 \not = 0$.  Then $b = {2a^2 \over 5}$,
$$c_2 = {3ab-15c \over 25} = {6a^3 - 75c \over 125},$$
$$g_2(t) = c_2t^2 + d_2t +e_2,$$
$$g_1(t) = 5t^4 + 4at^3 + {6a^2 \over 5}t^2 + 2ct + d.$$
$$g_1(t) - {5t^2 \over c_2}g_2(t) = \left(4a - {5d_2 \over c_2}\right)t^3 +\left({6a^2 \over 5}-{5e_2 \over c_2}\right)t^2 + 2ct + d,$$
\begin{eqnarray*}
&& g_1(t) - {5t^2 \over c_2}g_2(t) - \left(4a - {5d_2 \over c_2}\right){t \over c_2}g_2(t) \\
&= & \left[\left({6a^2 \over 5}-{5e_2 \over c_2}\right)-\left(4a - {5d_2 \over c_2}\right){d_2 \over c_2}\right]t^2 + \left[2c-\left(4a - {5d_2 \over c_2}\right){e_2 \over c_2}\right]t + d \\
&= & \left[{6a^2 \over 5}-{5e_2+4ad_2 \over c_2} + {5d_2^2 \over c_2^2}\right]t^2 + \left[2c-{4ae_2 \over c_2} + {5d_2e_2 \over c_2^2}\right]t + d \\
&= & {6a^2c_2^2 -5(5e_2+4ad_2)c_2 + 25d_2^2 \over 5c_2^2}t^2 + {2cc_2^2-4ae_2c_2 + 5d_2e_2 \over c_2^2}t + d \\
\end{eqnarray*}
\begin{eqnarray*}
&& g_1(t) - {5t^2 \over c_2}g_2(t) - \left(4a - {5d_2 \over c_2}\right){t \over c_2}g_2(t) - {6a^2c_2^2 -5(5e_2+4ad_2)c_2 + 25d_2^2 \over 5c_2^3}g_2(t) \\
&= & \left[{2cc_2^2-4ae_2c_2 + 5d_2e_2 \over c_2^2}-{6a^2c_2^2d_2 -5(5e_2+4ad_2)c_2d_2 + 25d_2^3 \over 5c_2^3}\right]t\\
& + & d - {6a^2c_2^2e_2 -5(5e_2+4ad_2)c_2e_2 + 25d_2^2e_2 \over 5c_2^3}\\
&= & {10cc_2^3-20ae_2c_2^2 + 25c_2d_2e_2 -6a^2c_2^2d_2 +5(5e_2+4ad_2)c_2d_2 - 25d_2^3 \over 5c_2^3}t\\
& + & {5c_2^3d - 6a^2c_2^2e_2 +5(5e_2+4ad_2)c_2e_2 - 25d_2^2e_2 \over 5c_2^3}.
\end{eqnarray*}
Then
$$g_3(t) = d_3t + e_{3,1},$$
where
$$d_3 = {\bar d_3 \over 5c_2^4},\ e_{3,1} = {\bar e_{3, 1} \over 5c_2^4},$$
$$\bar d_3 = -10cc_2^4+20ae_2c_2^3 -25c_2^2d_2e_2 + 6a^2c_2^3d_2 -5(5e_2+4ad_2)c_2^2d_2 + 25c_2d_2^3,$$
$$\bar e_{3,1} = - 5c_2^4d + 6a^2c_2^3e_2 -5(5e_2+4ad_2)c_2^2e_2 + 25c_2d_2^2e_2.$$

There are three subcases:

(3A) $d_3 = e_{3, 1} = 0$, i.e., $\bar d_3 = \bar e_{3, 1}=0$.   Then $g_2$ is the GCD of $g$ and $g'$.
Let $u = {d_2 \over c_2}$ and $v = {e_2 \over c_2}$.  We may use Propositions \ref{p2}, \ref{p3} and Theorem \ref{tt1} to determine the situation of $g(t)$.

(3B) $d_3 = 0$ and $e_{3, 1} \not = 0$, i.e., $\bar d_3 = 0$ and $\bar e_{3, 1} \not = 0$.  Then $g_3(t) = e_{3,1}$, $S = \{ g_0, g_1, g_2, g_3 \}$,
$$S(\infty) = \left\{+\infty, +\infty, c_2, e_{3,1} \right\},$$
$$S(0) = \{e, d, e_2, e_{3,1} \}.$$
We may replace $e_{3, 1}$ by $\bar e_{3, 1}$ here.
Then $V(0) = V(\infty)$ if and only if (i) $V(0) = V(\infty) =0$; or (ii) $V(0) = V(\infty) = 1$; or (iii) $V(\infty) = 2$, which implies $V(0) =2$ by the early discussion.

We may replace $e_{3, 1}$ by $\bar e_{3,1}$ here.
Thus, in this case, $g$ is strictly copositive if and only if either (i) $\bar e_{3, 1} > 0$, $c_2 \ge 0$,  $d \ge 0$, $e_2 \ge 0$; or (ii) $\bar e_{3, 1} < 0$, and $d \ge 0$ if $e_2 > 0$; or (iii) $\bar e_{3, 1} > 0$, $c_2 < 0$. %, $\min \{ d, e_2 \} < 0$.

(3C) $d_3 \not = 0$, i.e., $\bar d_3 \not = 0$. Then
$$g_2(t) - {c_2 \over d_3}tg_3(t) = \left(d_2-{c_2 \over d_3}e_{3,1}\right)t + e_2 = {d_2d_3 -c_2e_{3,1}\over d_3}t + e_2.$$
$$g_2(t) - {c_2 \over d_3}tg_3(t) - {d_2d_3 -c_2e_{3,1}\over d_3^2}g_3(t) = e_2 - {d_2d_3 -c_2e_{3,1}\over d_3^2}e_{3,1},$$
\begin{eqnarray*}
g_4(t) = e_4 & = & {d_2d_3 -c_2e_{3,1}\over d_3^2}e_{3,1} - e_2 \\
& = & {d_2d_3e_{3,1} - c_2e_{3,1}^2 - e_2d_3^2 \over d_3^2}\\
& = & {\bar e_4 \over \bar d_3^2},
\end{eqnarray*}
where
$$\bar e_4 = d_2\bar d_3\bar e_{3,1} - c_2\bar e_{3,1}^2 - \bar e_2\bar d_3^2.$$

There are two further subcases.

(3Ca) $e_4 = 0$, i.e., $\bar e_4 = 0$.   Then
 $g_3(t)$ is the GCD of $g$ and $g'$, and $\alpha = - {\bar e_{3, 1} \over \bar d_3}$ is a double root of $g$.  We may use Proposition \ref{p1} and Theorem \ref{tt1} to determine the situation of $g(t)$.

(3Cb) $e_4 \not = 0$, i.e., $\bar e_4 \not = 0$.  We have
$$S = \{ g_0, g_1, g_2, g_3, g_4 \},$$
$$S(\infty) = \left\{+\infty, +\infty, c_2, d_3, e_4 \right\},$$
$$S(0) = \{e, d, e_2, e_{3, 1}, e_4 \}.$$
We may  replace $d_3$, $e_{3, 1}$ and $e_4$ by $\bar d_3$, $\bar e_{3, 1}$ and $\bar e_4$ here.
Then $V(0) = V(\infty)$ if and only if (i) $V(0) = V(\infty) = 0$; or (ii) $V(0) = V(\infty) = 1$; or (iii) $V(0) = V(\infty) = 2$; or (iv) $V(\infty) = 3$, which implies $V(0) = 3$ by the early note.

Thus, in this case, $g$ is strictly copositive if and only if either (i) $\bar e_4 > 0$, $c_2 \ge 0$, $\bar d_3 \ge 0$, $d \ge 0$, $e_2 \ge 0$, $\bar e_{3, 1} \ge 0$; or (ii) $\bar e_4 < 0$, $c_2 \ge 0$ if $d_3 > 0$, $d \ge 0$ if $\max \{ e_2, \bar e_{3, 1} \}> 0$, $e_2 \ge 0$ if $\bar e_{3, 1} > 0$; or (iii) $\bar e_4 > 0$, $\min\{ c_2, \bar d_3 \} < 0$, $\min \{ d, e_2, \bar e_{3, 1} \} < 0$, $d \ge 0$ if $e_2 > 0$; or (iv) $\bar e_4 < 0$, $\bar d_3 > 0$, $c_2 < 0$.%, $\min\{ d, e_2 \} < 0$, $\max \{ e_2, \bar e_{3,1}\} > 0$.

(4) $b_2 \not = 0$.  Then
$$g_2(t) = b_2t^3 + c_2t^2 + d_2t + e_2.$$
$$g_1(t) - {5 \over b_2}g_2(t) = \left(4a - {5c_2 \over b_2}\right)t^3 + \left(3b - {5d_2 \over b_2}\right)t^2 + \left(2c - {5e_2 \over b_2}\right)t + d,$$
\begin{eqnarray*}
& & g_1(t) - {5 \over b_2}g_2(t) - {1 \over b_2}\left(4a - {5c_2 \over b_2}\right)g(t) \\
& = & \left[\left(3b - {5d_2 \over b_2}\right)- {c_2 \over b_2}\left(4a - {5c_2 \over b_2}\right)\right]t^2 + \left[\left(2c - {5e_2 \over b_2}\right)- {d_2 \over b_2}\left(4a - {5c_2 \over b_2}\right)\right]t + d- {e_2 \over b_2}\left(4a - {5c_2 \over b_2}\right)\\
& = & {3bb_2 - (5d_2+4ac_2)b_2 + 5c_2e_2 \over b_2^2}t^2 + {2cb_2^2 - (5e_2+4ad_2)b_2 + 5c_2d_2 \over b_2^2}t + {db_2^2- 4ae_2b_2 + 5c_2e_2 \over b_2^2}.\\
\end{eqnarray*}
Then
$$g_3(t) = c_3t^2 + d_{3, 1}t + e_{3, 2},$$
where
$$c_3 = {\bar c_3 \over b_2^2},\ d_{3,1} = {\bar d_{3, 1} \over b_2^2},\ e_{3, 2} = {\bar e_{3, 2} \over b_2^2},$$
$$\bar c_3 = - 3bb_2 + (5d_2+4ac_2)b_2 - 5c_2e_2,$$
$$\bar d_{3, 1} = - 2cb_2^2 + (5e_2+4ad_2)b_2 - 5c_2d_2,$$
$$\bar e_{3, 2} = -db_2^2+ 4ae_2b_2 - 5c_2e_2.$$

We have four subcases.

(4A) $c_3 = d_{3,1} = e_{3, 2} = 0$, i.e., $\bar c_3 = \bar d_{3, 1} = \bar e_{3, 2} = 0$.  This implies that $g_2$ is the GCD of $g$ and $g'$.  Let $p = {c_2 \over b_2}$, $q = {d_2 \over b_2}$ and $r = {e_2 \over b_2}$.   Then Proposition \ref{p4} determines the situation of $g(t)$.

(4B) $c_3 = d_{3,1}= 0$ and $e_{3, 2} \not = 0$, i.e., $\bar c_3 = \bar d_{3,1}= 0$ and $\bar e_{3, 2} \not = 0$.  Then
$$g_3(t) = e_{3, 2},$$
$$S = \{ g_0, g_1, g_2, g_3 \},$$
$$S(\infty) = \left\{+\infty, +\infty, b_2, e_{3, 2} \right\},$$
$$S(0) = \{e, d, e_2, e_{3, 2} \}.$$
We may replace $e_{3, 2}$ by $\bar e_{3, 2}$ here.
Then, $g$ is strictly copositive if and only if either (i) $\bar e_{3, 2} > 0$, $b_2 \ge 0$, $d \ge 0$, $e_2 \ge 0$; or (ii) $\bar e_{3, 2} < 0$, $d \ge 0$ if $e_2 > 0$; or (iii) $\bar e_{3, 2} > 0$, $b_2 < 0$.%, $\min \{ d, e_2 \} < 0$.

(4C) $c_3 = 0$ and $d_{3,1} \not = 0$.   Then
$$g_3(t) = d_{3, 1}t + e_{3, 2},$$
$$g_4(t) \equiv e_{4, 1} = -g_2\left(-{e_{3, 2} \over d_{3, 1}}\right)= {\bar e_{4, 1} \over \bar d_{3,1}^4},$$
where
$$\bar e_{4,1} = b_2\bar d_{3,1}\bar e_{3,2}^3-c_2\bar d_{3,1}^2\bar e_{3, 2}^2+d_2\bar d_{3,1}^3\bar e_{3,2}-e_2\bar d_{3,1}^4.$$

Again, there are two further subcases.

(4Ca) $e_{4, 1} = 0$, i.e., $\bar e_{4, 1} = 0$. Then $g_3(t)$ is the GCD of $g$ and $g'$, and $g$ has a double nonzero root
$$\alpha = -{\bar e_{3, 2} \over \bar d_{3, 1}}.$$
We may use Proposition \ref{p1} and Theorem \ref{tt1} to determine the situation of $g(t)$.

(4Cb) $e_{4, 1} \not = 0$, i.e., $\bar e_{4, 1} \not = 0$.  We have
$$S = \{ g_0, g_1, g_2, g_3, g_4 \},$$
$$S(\infty) = \left\{+\infty, +\infty, b_2, d_{3, 1}, e_{4, 1} \right\},$$
$$S(0) = \{e, d, e_2, e_{3, 2}, e_{4, 1} \}.$$
Then $V(0) = V(\infty)$ if and only if (i) $V(0) = V(\infty) = 0$; or (ii) $V(0) = V(\infty) = 1$; or (iii) $V(0) = V(\infty) = 2$; or (iv) $V(\infty) = 3$, which implies that $V(0) = 3$ by the early note.

We may replace $d_{3,1}$, $e_{3, 2}$ and $e_{4, 1}$ by $\bar d_{3,1}$, $\bar e_{3, 2}$ and $\bar e_{4, 1}$
Thus, in this subcase, $g$ is strictly copositive if and only if either (i) $\bar e_{4,1} > 0$, $b_2 \ge 0$, $\bar d_{3,1} \ge 0$, $d \ge 0$, $e_2 \ge 0$, $\bar e_{3, 2} \ge 0$; or (ii) $\bar e_{4,1} < 0$, $b_2 \ge 0$ if $\bar d_{3,1} > 0$, $d \ge 0$ if $\max \{ e_2, \bar e_{3, 2} \}> 0$, $e_2 \ge 0$ if $\bar e_{3, 2} > 0$; or (iii) $\bar e_{4,1} > 0$, $\min\{ b_2, \bar d_{3,1} \} < 0$, $\min \{ d, e_2, \bar e_{3, 2} \} < 0$, $d \ge 0$ if $e_2 > 0$; or (iv) $\bar e_{4,1} < 0$, $\bar d_{3,1} > 0$, $b_2 < 0$.%, $\min\{ d, e_2 \} < 0$, $\max \{ e_2, \bar e_{3,2}\} > 0$.

(4D) $c_3 \not = 0$, i.e., $\bar c_3 \not = 0$.
Then
$$g_3(t) = c_3t^2 + d_{3, 1}t + e_{3, 2},$$
$$g_2(t) = b_2t^3 + c_2t^2 + d_2t + e_2,$$
$$g_2(t) - {b_2 \over c_3}tg_3(t) = \left(c_2 - {b_2d_{3, 1} \over c_3}\right)t^2 + \left(d_2 - {b_2e_{3, 2} \over c_3}\right)t + e_2,$$
\begin{eqnarray*}
& & g_2(t) - {b_2 \over c_3}t g_3(t) - {1 \over c_3}\left(c_2 - {b_2d_{3, 1} \over c_3}\right)g_3(t) \\
& = & \left[\left(d_2 - {b_2e_{3, 2} \over c_3}\right)- {d_{3, 1} \over c_3}\left(c_2 - {b_2d_{3, 1} \over c_3}\right)\right]t + e_2 - {e_{3, 2} \over c_3}\left(c_2 - {b_2d_{3, 1} \over c_3}\right)\\
& = & \left(d_2 - {b_2e_{3, 2} \over c_3}- {c_2d_{3, 1} \over c_3} + {b_2d_{3, 1}^2 \over c_3^2}\right)t + e_2 - {c_2e_{3, 2} \over c_3} + {b_2d_{3, 1}e_{3, 2} \over c_3^2}.
\end{eqnarray*}
Then
$$g_4(t) = d_4t + e_{4, 2},$$
where
$$d_4 = {\bar d_4 \over \bar c_3^2},\ e_{4,2} = {\bar e_{4, 2} \over \bar c_3^2},$$
$$\bar d_4 = -\bar c_3^2d_2 + b_2\bar c_3\bar e_{3, 2} + c_2\bar c_3\bar d_{3, 1} - b_2\bar d_{3, 1}^2,$$
$$\bar e_{4, 2} = -\bar c_3^2e_2 + c_2\bar c_3\bar e_{3, 2} - b_2\bar d_{3, 1}\bar e_{3, 2}.$$

There are three further subcases.

(4Da) $d_4 = e_{4, 2} = 0$, i.e., $\bar d_4 = \bar e_{4, 2} = 0$.   Then $g_3$ is the GCD of $g$ and $g'$.
Let $u = {\bar d_{3, 1} \over \bar c_3}$ and $v = {\bar e_{3, 2} \over \bar c_3}$.  We may use
Propositions \ref{p2}, \ref{p3} and Theorem \ref{tt1} to determine the situation of $g(t)$.

(4Db) $d_4 = 0$ and $e_{4, 2} \not = 0$, i.e., $\bar d_4 = 0$ and $\bar e_{4, 2} \not = 0$.   Then $g_4(t) = e_{4, 2}$,
$$S = \{ g_0, g_1, g_2, g_3, g_4 \},$$
$$S(\infty) = \{ +\infty, + \infty, b_2, c_3, e_{4, 2} \},$$
$$S(0) = \{ e, d, e_2, e_{3, 2}, e_{4, 2} \}.$$
Then $V(0) = V(\infty)$ if and only if (i) $V(0) = V(\infty) = 0$; or (ii) $V(0) = V(\infty) = 1$; or (iii) $V(0) = V(\infty) = 2$; or (iv) $V(\infty) = 3$, which implies $V(0) = 3$ by the early note.

We may replace $c_3$, $e_{3, 2}$ and $e_{4, 2}$ here by $\bar c_3$, $\bar e_{3, 2}$ and $\bar e_{4, 2}$.
Thus, in this case, $g$ is strictly copositive if and only if either (i) $\bar e_{4,2} > 0$, $b_2 \ge 0$, $\bar c_3 \ge 0$, $d \ge 0$, $e_2 \ge 0$, $\bar e_{3, 2} \ge 0$; or (ii) $\bar e_{4,2} < 0$, $b_2 \ge 0$ if $\bar c_3 > 0$, $d \ge 0$ if $\max \{ e_2, \bar e_{3, 2} \}> 0$, $e_2 \ge 0$ if $\bar e_{3, 2} > 0$; or (iii) $\bar e_{4,2} > 0$, $\min\{ b_2, \bar c_3 \} < 0$, $\min \{ d, e_2, \bar e_{3, 2} \} < 0$, $d \ge 0$ if $e_2 > 0$; or (iv) $\bar e_{4,2} < 0$, $\bar c_3 > 0$, $b_2 < 0$.%, $\min\{ d, e_2 \} < 0$, $\max \{ e_2, \bar e_{3,2}\} > 0$.

(4Dc) $d_4 \not = 0$, i.e., $\bar d_4 \not = 0$.   Then
$$g_3(t) = c_3t^2 +d_{3, 1}t + e_{3, 2},$$
$$g_4(t) = d_4t + e_{4, 2},$$
$$g_5(t) \equiv e_5 = -g_3\left(-{e_{4, 2} \over d_4}\right)= {\bar e_5 \over b_2^2\bar d_4^2},$$
where
$$\bar e_5 = -\bar c_3\bar e_{4,2}^2 + \bar d_{3,1}\bar e_{4,2}\bar d_4-\bar e_{3,2}\bar d_4^2.$$

There are two further subcases.

(4Dc1) $e_5 = 0$, i.e., $\bar e_5 = 0$. Then $g_4(t)$ is the GCD of $g$ and $g'$, and $g$ has a double nonzero root
$$\alpha = -{\bar e_{4, 2} \over \bar d_4}.$$
We may use Proposition \ref{t1} and Theorem \ref{tt1} to determine the situation of $g(t)$.

(4Dc2) $e_5 \not = 0$, i.e., $\bar e_5 \not = 0$.  We have
$$S = \{ g_0, g_1, g_2, g_3, g_4, g_5 \},$$
$$S(\infty) = \left\{+\infty, +\infty, b_2, c_3, d_4, e_5 \right\},$$
$$S(0) = \{e, d, e_2, e_{3, 2}, e_{4, 2}, e_5 \}.$$
Then $V(0) = V(\infty)$ if and only if (i) $V(0) = V(\infty) = 0$; or (ii) $V(0) = V(\infty) = 1$; or (iii) $V(0) = V(\infty) = 2$; or (iv) $V(0) = V(\infty) = 3$, or (v) $V(\infty) = 4$, which implies that $V(0) = 4$ by the early note.

We may replace $c_3, d_4, e_{3, 2}, e_{4, 2}$ and $e_5$ here by $\bar c_3, \bar d_4, \bar e_{3, 2}, \bar e_{4, 2}$ and $\bar e_5$.
Thus, in this case, $g$ is strictly copositive if and only if either (i) $\bar e_5 > 0$, $b_2 \ge 0$, $\bar c_3 \ge 0$, $\bar d_4 \ge 0$, $\bar e_{4, 2} \ge 0$, $d \ge 0$, $e_2 \ge 0$, $\bar e_{3, 2} \ge 0$; or (ii) $\bar e_5< 0$, $b_2 \ge 0$ if $\max\{ \bar c_3, \bar d_4\} > 0$, $\bar c_3 \ge 0$ if $\bar d_4 > 0$, $d \ge 0$ if $\max \{ e_2, \bar e_{3, 2}, \bar e_{4, 2} \} > 0$, $e_2 \ge 0$ if $\max \{ \bar e_{3,2}, \bar e_{4,2} \} > 0$, $\bar e_{3, 2} \ge 0$ if $\bar e_{4,2} > 0$; or (iii) $\bar e_5 > 0$, $\min \{ b_2, \bar c_3, \bar d_4 \} < 0$, $b_2 \ge 0$ if $\bar c_3 > 0$, $\min \{d, e_2, \bar e_{3, 2}, \bar e_{4,2} \} < 0$, $d \ge 0$ if $\max \{ e_2, e_{3,2} \} > 0$, $e_2 \ge 0$ if $e_{3, 2} > 0$; or (iv) $e_5 < 0$, $\min \{ b_2, \bar c_3 \} < 0$, $\max \{\bar c_3, \bar d_4 \} > 0$, $\min \{ d, e_2, \bar e_{3,2} \} < 0$, $\max \{ e_2, \bar e_{3,2}, \bar e_{4,2} \} > 0$, $d \ge 0$ if $e_2 > 0$; or (v) $\bar e_5 > 0$, $\bar d_4 < 0$, $\bar c_3 > 0$, $b_2 < 0$.%, $\min\{\bar e_{3,2}, \bar e_{4,2}\} <0$, $\max \{e_2, \bar e_{3,2} \} > 0$, $\min\{d, e_2 \} < 0$.

\subsection{Positivity Conditions for Quintic Polynomials}

Summarizing the coefficients, we have
\begin{equation} \label{e2}
b_2 = {4a^2-10b \over 25},\ c_2 = {3ab-15c \over 25}, d_2 = {2ac-20d \over 25}, e_2 = {ad-25e \over 25},
\end{equation}
%Under the condition that $b_2 = c_2 = 0$ and $d_2 \not = 0$, we have
\begin{equation} \label{e3}
\bar e_3 = -125e_2^4+100ae_2^3d_2-30a^2e_2^2d_2^2+4a^3e_2d_2^3-dd_2^4,
\end{equation}
%Under the condition that $b_2 = 0$ and $c_2 \not = 0$, we have
\begin{equation} \label{e4}
\bar d_3 = -10cc_2^4+20ae_2c_2^3 -25c_2^2d_2e_2 + 6a^2c_2^3d_2 -5(5e_2+4ad_2)c_2^2d_2 + 25c_2d_2^3,
\end{equation}
\begin{equation} \label{e5}
\bar e_{3,1} = - 5c_2^4d + 6a^2c_2^3e_2 -5(5e_2+4ad_2)c_2^2e_2 + 25c_2d_2^2e_2,
\end{equation}
%Under the further condition that $\bar d_3 \not = 0$, we have
\begin{equation} \label{e6}
\bar e_4 = d_2\bar d_3\bar e_{3,1} - c_2\bar e_{3,1}^2 - \bar e_2\bar d_3^2,
\end{equation}
%Under the condition $b_2 \not = 0$, we have
\begin{equation} \label{e7}
\bar c_3 = -3bb_2 + (5d_2+4ac_2)b_2 - 5c_2e_2,
\end{equation}
\begin{equation} \label{e8}
\bar d_{3, 1} = -2cb_2^2 +(5e_2+4ad_2)b_2 - 5c_2d_2,
\end{equation}
\begin{equation} \label{e9}
\bar e_{3, 2} = -db_2^2+ 4ae_2b_2 - 5c_2e_2,
\end{equation}
%Under the further condition that $\bar c_3 = 0$ and $\bar d_{3, 1} \not = 0$, we have
\begin{equation} \label{e10}
\bar e_{4,1} = b_2\bar d_{3,1}\bar e_{3,2}^3-c_2\bar d_{3,1}^2\bar e_{3, 2}^2+d_2\bar d_{3,1}^3\bar e_{3,2}-e_2\bar d_{3,1}^4,
\end{equation}
%Under the further condition that $\bar c_3 \not = 0$, we have
\begin{equation} \label{e11}
\bar d_4 = -\bar c_3^2d_2 + b_2\bar c_3\bar e_{3, 2} + c_2\bar c_3\bar d_{3, 1} - b_2\bar d_{3, 1}^2,
\end{equation}
\begin{equation} \label{e12}
\bar e_{4, 2} = -\bar c_3^2e_2 + c_2\bar c_3\bar e_{3, 2} - b_2\bar d_{3, 1}\bar e_{3, 2},
\end{equation}
%Under the further condition that $\bar c_3 \not = 0$ and $\bar d_4 \not = 0$, we have
\begin{equation} \label{e13}
\bar e_5 = -\bar c_3\bar e_{4,2}^2 + \bar d_{3,1}\bar e_{4,2}\bar d_4-\bar e_{3,2}\bar d_4^2.
\end{equation}

We now can state our theorem.

\begin{theorem} \label{t1}
Let $g(t)$ be defined by (\ref{e1}), with $e > 0$. Let the additional coefficients be defined by (\ref{e2}-\ref{e13}).  Then we have $11$ cases.
%$g(t) > 0$ for all $t \ge 0$ if and only if one of the following $15$ conditions holds.

(1) $b_2 = c_2 = d_2= 0$. In this case, $g(t) > 0$ for all $t \ge 0$.

(2)  $b_2 = c_2 = 0$, $d_2 \not = 0$ and $\bar e_3 \not = 0$.  In this case, $g(t) > 0$ for all $t \ge 0$ if and only if
either (i) $\bar e_3 > 0$, $d_2 > 0$, $d \ge 0$, $e_2 \ge 0$; or (ii) $\bar e_3 < 0$, and $d \ge 0$ if $e_2 > 0$; or (iii) $\bar e_3 > 0$, $d_2 < 0$. %,  $\min \{ d, e_2 \} < 0$.

(3) $b_2 = 0$, $c_2 \not = 0$, $\bar d_3 = 0$ and $\bar e_{3, 1} \not = 0$.  In this case, $g(t) > 0$ for all $t \ge 0$ if and only if either (i) $\bar e_{3, 1} > 0$, $c_2 \ge 0$,  $d \ge 0$, $e_2 \ge 0$; or (ii) $\bar e_{3, 1} < 0$, and $d \ge 0$ if $e_2 > 0$; or (iii) $\bar e_{3, 1} > 0$, $c_2 < 0$.%, $\min \{ d, e_2 \} < 0$.

(4) $b_2 = 0$, $c_2 \not = 0$, $\bar d_3 \not = 0$ and $\bar e_4 \not = 0$.  In this case, $g(t) > 0$ for all $t \ge 0$ if and only if either (i) $\bar e_4 > 0$, $c_2 \ge 0$, $\bar d_3 \ge 0$, $d \ge 0$, $e_2 \ge 0$, $\bar e_{3, 1} \ge 0$; or (ii) $\bar e_4 < 0$, $c_2 \ge 0$ if $\bar d_3 > 0$, $d \ge 0$ if $\max \{ e_2, \bar e_{3, 1} \}> 0$, $e_2 \ge 0$ if $\bar e_{3, 1} > 0$; or (iii) $\bar e_4 > 0$, $\min\{ c_2, \bar d_3 \} < 0$, $\min \{ d, e_2, \bar e_{3, 1} \} < 0$, $d \ge 0$ if $e_2 > 0$; or (iv) $\bar e_4 < 0$, $\bar d_3 > 0$, $c_2 < 0$.%, $\min\{ d, e_2 \} < 0$, $\max \{ e_2, \bar e_{3,1}\} > 0$.

(5) $b_2 \not = 0$, $\bar c_3 = \bar d_{3,1}= 0$, and $\bar e_{3, 2} \not = 0$.  In this case, $g(t) > 0$ for all $t \ge 0$ if and only if either (i) $\bar e_{3, 2} > 0$, $b_2 \ge 0$, $d \ge 0$, $e_2 \ge 0$; or (ii) $\bar e_{3, 2} < 0$, $d \ge 0$ if $e_2 > 0$; or (iii) $\bar e_{3, 2} > 0$, $b_2 < 0$.%, $\min \{ d, e_2 \} < 0$.

(6) $b_2 \not = 0$, $\bar c_3 = 0$, $\bar d_{3,1} \not = 0$ and $\bar e_{4,1} \not = 0$.   In this case, $g(t) > 0$ for all $t \ge 0$ if and only if either (i) $\bar e_{4,1} > 0$, $b_2 \ge 0$, $\bar d_{3,1} \ge 0$, $d \ge 0$, $e_2 \ge 0$, $e_{3, 2} \ge 0$; or (ii) $\bar e_{4,1} < 0$, $b_2 \ge 0$ if $\bar d_{3,1} > 0$, $d \ge 0$ if $\max \{ e_2, \bar e_{3, 2} \}> 0$, $e_2 \ge 0$ if $\bar e_{3, 2} > 0$; or (iii) $\bar e_{4,1} > 0$, $\min\{ b_2, \bar d_{3,1} \} < 0$, $\min \{ d, e_2, \bar e_{3, 2} \} < 0$, $d \ge 0$ if $e_2 > 0$; or (iv) $\bar e_{4,1} < 0$, $\bar d_{3,1} > 0$, $b_2 < 0$.%, $\min\{ d, e_2 \} < 0$, $\max \{ e_2, \bar e_{3,2}\} > 0$.

(7) $b_2 \not = 0$, $\bar c_3 \not = 0$ and $\bar d_4 = 0$.  In this case, $g(t) > 0$ for all $t \ge 0$ if and only if either (i) $\bar e_{4,2} > 0$, $b_2 \ge 0$, $\bar c_3 \ge 0$, $d \ge 0$, $e_2 \ge 0$, $\bar e_{3, 2} \ge 0$; or (ii) $\bar e_{4,2} < 0$, $b_2 \ge 0$ if $\bar c_3 > 0$, $d \ge 0$ if $\max \{ e_2, \bar e_{3, 2} \}> 0$, $e_2 \ge 0$ if $\bar e_{3, 2} > 0$; or (iii) $\bar e_{4,2} > 0$, $\min\{ b_2, \bar c_3 \} < 0$, $\min \{ d, e_2, \bar e_{3, 2} \} < 0$, $d \ge 0$ if $e_2 > 0$; or (iv) $\bar e_{4,2} < 0$, $\bar c_3 > 0$, $b_2 < 0$.%, $\min\{ d, e_2 \} < 0$, $\max \{ e_2, \bar e_{3,2}\} > 0$.

(8) $b_2 \not = 0$, $\bar c_3 \not = 0$, $\bar d_4 \not = 0$ and $\bar e_5 \not = 0$.  In this case, $g(t) > 0$ for all $t \ge 0$ if and only if  either (i) $\bar e_5 > 0$, $b_2 \ge 0$, $\bar c_3 \ge 0$, $\bar d_4 \ge 0$, $\bar e_{4, 2} \ge 0$, $d \ge 0$, $e_2 \ge 0$, $\bar e_{3, 2} \ge 0$; or (ii) $\bar e_5< 0$, $b_2 \ge 0$ if $\max\{ \bar c_3, \bar d_4\} > 0$, $\bar c_3 \ge 0$ if $\bar d_4 > 0$, $d \ge 0$ if $\max \{ e_2, \bar e_{3, 2}, \bar e_{4, 2} \} > 0$, $e_2 \ge 0$ if $\max \{ \bar e_{3,2}, \bar e_{4,2} \} > 0$, $\bar e_{3, 2} \ge 0$ if $\bar e_{4,2} > 0$; or (iii) $\bar e_5 > 0$, $\min \{ b_2, \bar c_3, \bar d_4 \} < 0$, $b_2 \ge 0$ if $\bar c_3 > 0$, $\min \{d, e_2, \bar e_{3, 2}, \bar e_{4,2} \} < 0$, $d \ge 0$ if $\max \{ e_2, \bar e_{3,2} \} > 0$, $e_2 \ge 0$ if $\bar e_{3, 2} > 0$; or (iv) $\bar e_5 < 0$, $\min \{ b_2, \bar c_3 \} < 0$, $\max \{\bar c_3, \bar d_4 \} > 0$, $\min \{ d, e_2, \bar e_{3,2} \} < 0$, $\max \{ e_2, \bar e_{3,2}, \bar e_{4,2} \} > 0$, $d \ge 0$ if $e_2 > 0$; or (v) $\bar e_5 > 0$, $\bar d_4 < 0$, $\bar c_3 > 0$, $b_2 < 0$.%, $\min\{\bar e_{3,2}, \bar e_{4,2}\} <0$, $\max \{e_2, \bar e_{3,2} \} > 0$, $\min\{d, e_2 \} < 0$.

In these eight cases, $g(t) > 0$ for all $t \ge 0$ if and only if $g(t) \ge 0$ for all $t \ge 0$.

(9) (i) $b_2 = c_2 = 0$, $d_2 \not = 0$ and $\bar e_3 = 0$. Let $\alpha = -{e_2 \over d_2}$.

(ii) $b_2 = 0$, $c_2 \not = 0$, $\bar d_3 \not = 0$ and $\bar e_4 = 0$.  Let $\alpha = - {\bar e_{3, 1} \over \bar d_3}$.

(iii) $b_2 \not = 0$, $\bar c_3 = 0$, $\bar d_{3, 1} \not = 0$ and $\bar e_{4,1} = 0$.  Let $\alpha = -{\bar e_{3, 2} \over \bar d_{3, 1}}$.

(iv) $b_2 \not = 0$, $\bar c_3 \not = 0$, $\bar d_4 \not = 0$ and $\bar e_5 = 0$.  Let $\alpha = -{\bar e_{4, 2} \over \bar d_4}$.

Then we may apply Proposition \ref{p1} and Theorem \ref{tt1} to determine the situation of $g(t)$.

(10) (i) $b_2 = 0$, $c_2 \not = 0$ and $\bar d_3 = \bar e_{3, 1} = 0$.  Let $u = {d_2 \over c_2}$ and $v = {e_2 \over c_2}$.

(ii) $b_2 \not = 0$, $\bar c_3 \not = 0$ and $\bar d_4 = \bar e_{4, 2} = 0$.   Let $u = {\bar d_{3, 1} \over \bar c_3}$ and $v = {\bar e_{3, 2} \over \bar c_3}$.

Then we may use Propositions \ref{p2}, \ref{p3} and Theorem \ref{tt1} to determine the situation of $g(t)$.

(11) $b_2 \not = 0$ and $\bar c_3 = \bar d_{3,1} = \bar e_{3, 2} = 0$.  Let $p = {c_2 \over b_2}$, $q = {d_2 \over b_2}$ and $r = {e_2 \over b_2}$.   Then Proposition \ref{p4} determines the situation of $g(t)$.
\end{theorem}
{\bf Proof} The $11$ cases are summarized from the discussion in the last subsection. Cases (1-8) are corresponding to Cases (1), (2B), (3B), (3Cb), (4B), (4Cb), (4Db) and (4Dc2) of the last subsection, respectively.   Case (9) is corresponding Cases (2A), (3Ca), (4Ca) and (4Dc1) of the last subsection.   Case (10) is corresponding Cases (3A) and (4Da) of the last subsection.   Case (11) is corresponding to Case (4A) of the last subsection.
\qed

\section{Higher Even Degree Polynomials and Appendices}

We consider non-degenerate even degree polynomials with their degrees higher than or equal to $4$.
Suppose that
$$g(t) = \sum_{i=1}^m a_it^{m-i},$$
where $m \ge 4$ is even, $a_0 = a_m = 1$, $a_1, \cdots, a_{m-1}$ are real numbers.  Denote
$$\S = \{ (a_1, \cdots, a_{m-1}) \in \Re^{m-1} : g(t) > 0\ {\rm for\ all}\ t \},$$
$$\bar \S = \{ (a_1, \cdots, a_{m-1}) \in \Re^{m-1} : g(t) \ge 0\ {\rm for\ all}\ t \}$$
and
$$A = \{ (a_1, \cdots, a_{m-1}) \in \Re^{m-1} : g(t) = (t^2 + ut + v)^2\phi(t), u^2 < 4v, \phi(t) \ge 0\ {\rm for\ all}\ t \}.$$
We call $A$ the {\bf appendix} of $\bar S$.  We have the following theorem.

\begin{theorem} \label{at1}
In the above setting, $\bar S$ is a closed convex set, $S$ is its interior, and $A \subset S$.   The dimension of $\bar S$ and $S$ is $m-1$, while the dimension of $A$ is $m-3$.  The discriminant $\Delta(g)$ is equal to zero on $A$ and the boundary of $\bar S$, and has the same sign at the other part of $\bar S$.
\end{theorem}
{\bf Proof} Suppose that
$$g(t) = \sum_{i=1}^m a_it^{m-i},\ \hat g(t) = \sum_{i=1}^m \hat a_it^{m-i},$$
$g(t), \hat g(t) \ge 0$ for all $t$, $a_0 = a_m = \hat a_0 = \hat a_m =1$.   Then ${1 \over 2}(g(t) + \hat g(t)) \ge 0$ for all $t$.   This shows that $\bar S$ is convex.  Taking limiting points, we see that $\bar S$ is closed.   Similarly, $S$ is also convex.   Consider $\min \{ g(t) \}$.   By the continuity property, we see that $S$ is an open set, and $\bar S$ is the closure of $S$.  Then $S$ is the interior of $\bar S$.
For any $(a_1, \cdots, a_{m-1}) \in A$, $g(t) > 0$ for all $t$.  Thus, $A \subset S$.
We also see that there is an $\epsilon > 0$ such that $(a_1, \cdots, a_{m-1}) \in S$ for $|a_i| \le \epsilon, i = 1, \cdots, m-1$.   This shows that the dimension of $\bar S$ and $S$ is $m-1$.  Consider the number of independent parameters of $g(t)$ in $A$, we conclude that the dimension of $A$ is $m-3$.    For
$g(t)$ in $A$, $g(t)$ has multiple roots.  Thus $\Delta(g) = 0$.   On a boundary point of $\bar S$, as it neighbors some parts not in $\bar S$, we also have $\Delta(g) = 0$.  On the other part of $\bar S$, as it is in the interior of $\bar S$, and $g(t)$ has no multiple roots there, we have $\Delta(g) \not =0$.   As $A$ is of dimension $m-3$, the other parts is connected.  Thus, $\Delta(g)$ has the same sign there.
\qed

This theorem shows that, in the above setting, for all $g(t)$ without real roots, $\Delta(g)$ takes the same sign, which depends upon $m$ only, except at an appendix set of dimension lower by two, where $\Delta(g) = 0$.

Note that such an appendix also exists for odd degree polynomials with their degrees higher than $4$.
In Subsection 5.1, for the subcase (B1), if $u^2 < 4v$, then there exists also an appendix similarly.
However, as the property that $g(t) \ge 0$ for all $t$ only works for even degree polynomials, only for even degree polynomials, the properties of the appendix are distinguished.

\bigskip

{\bf Acknowledgment}  We are thankful to Chen Ouyang and Jinjie Liu for drawing Figures 1 and 2.


\bigskip


\begin{thebibliography}{99}

\bibitem{BPM06} S. Basu, R. Pollack and R. Marie-Fran\c{c}ois, {\sl Algorithm in Real Algebraic Geometry}, 2nd ed., Springer, Berlin, 2006.

\bibitem{CHQ17} H. Chen, Z. Huang and L. Qi, ``Copositivity detection of tensors: Theory and algorithm'', {\sl J. Optim. Theory Appl. \bf 174} (2017) 746-761.

\bibitem{CHQ18} H. Chen, Z. Huang and L. Qi, ``Copositive tensor detection and its application in physics and hypergraphs'', {\sl Comput. Optim. Appl. \bf 69} (2018) 133-158.

\bibitem{CW18} H. Chen and Y. Wang, ``Higher order copositive tensors and its application'', {\sl J. Appl. Anal. Comput. \bf 8} (2018) 1863-1885.

%\bibitem{DQW17} W. Ding, L. Qi and Y. Wei, ``Inheritance properties and sum-of-squares decomposition of Hankel tensors: Theory and algorithms'', {\sl BIT \bf 57} (2017) 169-190.

\bibitem{FI19} F.S. Faro and I.P. Ivanov, ``Boundedness from below in the $U(1) \times U(1)$ three-Higgs-doublet model'', {\sl Phys. Rev. D \bf 100} (2019) 035038.

\bibitem{GKZ94} I.M. Gelfand, M.M. Kapranov and A.V. Zelevinsky, {\sl Discriminants, Resultants and Multidimensional Determinants}, Birkh\"{a}user, Boston, 1994.

\bibitem{GL64} R.N. Gadem and C.C. Li, ``On positive definiteness of quartic forms of two variables'', {\sl IEEE Trans. Automat. Control \bf AC-9} (1964) 187-188.

\bibitem{Ga20} Y. Gao, ``A necessary and sufficient condition for the positive definite problem of a binary quartic form'', arXiv:2009.01033, 2020.

%\bibitem{HQ19} Z. Huang and L. Qi, ``Tensor complementarity problems - Part I: Basic Theory'', {\sl J. Optim. Theory Appl. \bf 183} (2019) 1-23.

%\bibitem{HQ19a} Z. Huang and L. Qi, ``Tensor complementarity problems - Part III: Applications'', {\sl J. Optim. Theory Appl. \bf 183} (2019) 771-791.

\bibitem{IKM18} I.P. Ivanov, M. K\"{o}pke and M. M\"{u}hlleitner, ``Boundedness from below in the $U(1) \times U(1)$ three-Higgs-doublet model'', {\sl Phys. Rev. D \bf 100} (2018) 035038.

\bibitem{JM81} E.I. Jury and M. Mansour, ``Positivity and nonnegativity of a quartic equation and related problems'', {\sl IEEE Trans. Automat. Control \bf 26} (1981) 444-451.

\bibitem{Ka16} K. Kannike, ``Vacuum stability of a general scalar potential of a few fields'', {\sl Eur. Phys. J. C \bf 76} (2016) 324.

\bibitem{Ka18} K. Kannike, ``Erratum to: Vacuum stability of a general scalar potential of a few fields'', {\sl Eur. Phys. J. C \bf 78} (2018) 355.

\bibitem{Ku65} W.H. Ku, ``Explicit criterion for the positive definiteness of a general quartic form'', {\sl IEEE Tram. Automat. Control \bf AC-10} (1965) 372-373.

\bibitem{LZHQ19} L. Li, X. Zhang, Z. Huang and L. Qi, ``Test of copositive tensors'', {\sl J. Indust. Manag. Optim. \bf 15} (2019) 881-891.

\bibitem{LS19} J. Liu and Y. Song, ``Analytical expressions of copositivity for 3rd-order symmetric tensors and applications'', arXiv:1911.10284, (2019).

%\bibitem{LQ16} Z. Luo and L. Qi, ``Completely positive tensors: Properties, easily checkable subclasses, and tractable relaxations'', {\sl SIAM J. Matrix Anal. Appl. \bf 37} （2016） 1675-1698.

\bibitem{NYZ18} J. Nie, Z. Yang,  and X. Zhang, ``A complete semi-definite algorithm for detecting copositive matrices and tensors'', {\sl SIAM J. Optim. \bf 28} (2018) 2902-2921.

%\bibitem{Qi05} L. Qi, ``Eigenvalues of a real supersymmetric tensor'', {\sl J. Symbolic Comput. \bf 40} (2005) 1302-1324.

\bibitem{Qi13} L. Qi, ``Symmetric nonnegative tensors and copositive tensors'', {\sl Linear Algebra Appl. \bf 439} (2013) 228-238.

%\bibitem{Qi15} L. Qi, ``Hankel tensors: Associated Hankel matrices and Vandermonde decomposition'', {\sl Commun. Math. Sci. \bf 13} (2015) 113-125.

\bibitem{QCC18} L. Qi, H. Chen and Y. Chen, {\sl Tensor Eigenvalues and Their Applications}, Springer, New York, 2018.

%\bibitem{QH19} L. Qi and Z. Huang, ``Tensor complementarity problems - Part II: Solution methods'', {\sl J. Optim. Theory Appl. \bf 183} (2019) 365-385.

\bibitem{QL17} L. Qi and Z. Luo, {\sl Tensor Analysis: Spectral Theory and Special Tensors}, SIAM, Philadelphia, 2017.

%\bibitem{QS14} L. Qi and Y. Song, ``An even order symmetric B tensor is positive defnite'', {\sl Linear Algebra Appl. \bf 457} (2014) 303-312.

%\bibitem{QXX14} L. Qi, C. Xu and Y. Xu, ``Nonnegative tensor factorization, completely positive tensors, and a hierarchical elimination algorithm'', {\sl SIAM J. Matrix Anal. Appl. \bf 35} （2014） 1227-1241.

\bibitem{SH88} J.W. Schmidt and W. He\ss, ``Positivity of cubic polynomials on intervals and positive spline interpolation'', {\sl BIT Numer. Math. \bf 28} (1988) 340-352.


\bibitem{SQ15} Y. Song and L. Qi, ``Necessary and sufficient conditions of copositive tensors'', {\sl Linear Multilinear Algebra \bf 63} (2015) 120-131.

%\bibitem{SQ16} Y. Song and L. Qi, ``Eigenvalue analysis of constrained minimization problems for homogeneous polynomial'', {\sl J. Global Optim. \bf 64} (2016) 563-575.

\bibitem{SQ20} Y. Song and L. Qi, ``Analytical expressions of copositivity for fourth-order symmetric tensors'', {\sl Anal. Appl.} DOI:10.1142/S0219530520500049 (2020).

\bibitem{UW94} G. Ulrich and L.T. Watson, ``Positivity conditions for quartic polynomials'', {\sl SIAM J. Sci. Comput. \bf 15} (1994) 528-544.

\bibitem{WQ05} F. Wang and L. Qi, ``Comments on `Explicit criterion for the positive definiteness of a general quartic form''', {\sl IEEE Tram. Automat. Control \bf 50} (2005) 416-418.



































%%%%%%%%%%%%%%%%%%%%%%%%%%%%%%%%%%%%%%%%%%%%%%%%%%%%%%%%%%%%%%%%%%%%%%%%

%

%

%

%

%

%

%

%

%

%

%

%

%

%

%

%

%

%

%

%

%

%

%

%

%

\end{thebibliography}
\end{document}